\documentclass[12pt]{article}

\usepackage{amsmath,amssymb,theorem,makeidx,latexsym,epsfig,url}

\newtheorem{defn}{Definition}[section]
\newtheorem{lemma}[defn]{Lemma}

{\theorembodyfont{\rmfamily}
\newtheorem{ex}[defn]{Example}}
\newtheorem{thm}[defn]{Theorem}
\newtheorem{prop}[defn]{Proposition}
\newtheorem{cor}[defn]{Corollary}
\newtheorem{rem}[defn]{Remark}

\def\bp{{\noindent\bf Proof. \ }}
\def\ep{\hfill$\square$\par\bigskip}

\newcommand{\h}{{\cal H}}
\newcommand{\lt}{{\ell^2}}
\newcommand{\mr}{\mathbb R}
\def\bt{\begin{thm}}
\def\et{\end{thm}}
\def\bc{\begin{cor}}
\def\ec{\end{cor}}
\def\bl{\begin{lemma}}
\def\el{\end{lemma}}

\newcommand{\ltr}{ L^2(\mathbb R) }

\def\la{\langle}
\def\ra{\rangle}
\def\ga{\gamma}
\def\ftk{\{f_k\}}
\def\htk{\{h_k\}}
\def\ctk{\{c_k\}}

\def\gtk{\left\{g_k\right\}}

\newcommand{\afh}{ \forall f \in \h}
\newcommand{\mn}{\mathbb N}

\newcommand{\mz}{\mathbb Z}
\newcommand{\mc}{\mathbb C}

\def\mts{\{E_{mb}T_{na}g\}_{m,n\in \mz}}
\def\bee{\begin{eqnarray}}
\def\ene{\end{eqnarray}}

\def\bes{\begin{eqnarray*}}
\def\ens{\end{eqnarray*}}

\def\bex{\begin{ex} }
\def\enx{\end{ex}}
\def\bei{\begin{itemize}}
\def\eni{\end{itemize}}

\title{Approximately dual frame pairs in Hilbert spaces and applications
to Gabor frames}

\date{\today}

\author{Ole Christensen and Richard S. Laugesen}

\begin{document}

\maketitle

\renewcommand{\thefootnote}{}
\footnote{\noindent \hspace*{-18pt}{\bf AMS MSC:} 42C15. {\bf Key
Words:} \text{Frames, approximate duals, Gabor frames, Gaussian.}}

\vspace*{-24pt}
\begin{abstract} We discuss the concepts of
pseudo-dual frames and approximately dual frames, and illuminate
their relationship to classical frames. Approximately dual frames
are easier to construct than the classical dual frames, and might
be tailored to yield almost perfect reconstruction.

For approximately dual frames constructed via perturbation theory,
we provide a bound on the deviation from perfect reconstruction. An
alternative bound is derived for the rich class of Gabor frames, by
using the Walnut representation of the frame operator to estimate
the deviation from equality in the duality conditions.

As illustration of the results, we construct explicit approximate
duals of Gabor frames generated by the Gaussian; these approximate
duals yield almost perfect reconstruction. Amazingly, the method
applies also to certain Gabor frames that are far from being tight.
\end{abstract}

\section{Introduction} \label{sec1}
Let $\h$ be a separable Hilbert space with inner product $\langle
\cdot , \cdot \rangle$. Given a frame $\ftk$ for $\h,$ it is well
known that there exists at least one dual frame $\gtk,$ \emph{i.e.},
a frame for which \bes f= \sum \la f,f_k\ra g_k, \quad \forall f \in
\h.\ens  Unfortunately, it is usually complicated to calculate a
dual frame explicitly. Hence we seek methods for constructing {\it
approximate} duals. The current paper treats the concepts of
pseudo-dual and approximately dual frames and examines their
properties. The idea of looking at approximately dual frames has
appeared several times in the literature, for example
\cite{BL,GHHLWW,H} in the wavelet setting and \cite{FK} for Gabor
systems, but our contribution is the first systematic treatment.

An approximately dual frame $\gtk$ associated with $\ftk$ satisfies
\bee \label{10a} \lVert f- \sum \la f,f_k\ra g_k \rVert \le \epsilon
\lVert f  \rVert, \quad \forall f\in \h,\ene for some $\epsilon <1$.
(This bound is particularly interesting when $\epsilon$ is small.)
We will show that any approximately dual frame generates a natural
dual frame in the classical sense; further we obtain a family of
frames that interpolates between the approximately dual frame and
the ``true" dual frame.

We also investigate the use of perturbation theory to construct
approximately dual frames. There are situations where it is hard to
find a dual frame for a given frame $\ftk,$ whereas a frame $\htk$
close to $\ftk$ can be found for which a dual frame $\gtk$ is known
explicitly. We present conditions under which such a frame $\gtk$ is
approximately dual to $\ftk$, both in the case where $\gtk$ is an
arbitrary dual of $\htk$ and for the particular case where it is the
canonical dual.

General Hilbert space estimates might of course be improved in
concrete cases. For Gabor frames in $\ltr$ we present a direct way
to obtain an inequality of the type \eqref{10a}, based on the Walnut
representation of the frame operator \cite{Wa}; the obtained value
$\epsilon$ measures the deviation from equality in the Gabor duality
conditions.

As an illustration of the method, we construct approximately dual
frames for two Gabor frames generated by the Gaussian. In these
concrete cases, we obtain approximations of type \eqref{10a} with
very small values for $\epsilon$.

The paper is organized as follows. Section \ref{sec2} contains a few
elementary definitions and results from standard frame theory. In
Section \ref{sec3} we introduce the new concepts of pseudo-dual and
approximately dual frames, and investigate their properties. The
(apparently) most useful concept of these, approximately dual
frames, is studied in more detail in Section \ref{sec4}; in
particular, we show how to obtain approximate duals via perturbation
theory. The direct estimates for Gabor frames are found in Section
\ref{gab}, and the applications to Gabor frames generated by the
Gaussian are in Section \ref{sec16}. Finally, an appendix collects a
few examples and proofs related to classical frame theory.

\section{Basic frame theory} \label{sec2}
For the introductory frame material that follows, see any standard
reference on frames, such as \cite{Da2}, \cite{HW}, \cite{Y},
\cite[Chapter~5]{CBN}.

A sequence $\ftk$ in $\h,$ indexed by an arbitrary countable index
set, is a {\it frame} if there exist constants $A,B>0$ such that
\bee \label{8a} A \lVert f \rVert^2 \le \sum | \langle f , f_k
\rangle |^2 \leq B \lVert f \rVert^2 \qquad \forall f \in \h. \ene
Any  numbers $A,B$ such that \eqref{8a} holds are called {\it (lower
and upper) frame bounds}.

If at least the upper inequality holds, then $\ftk$ is a {\it
Bessel sequence.} For later use, note that if $\{ f_k \}$ is a
Bessel sequence then so is $\{ W f_k \}$ for any bounded operator
$W : \h \to \h$.

Given a Bessel sequence $\ftk,$ the \emph{synthesis operator} $T :
\lt \to \h$ given by
\[
T \{ c_k \} = \sum c_k f_k\] is linear and bounded, with $\lVert T
\rVert \le \sqrt{B}$; the series converges unconditionally for all
$\ctk \in \ell^2.$ The adjoint of $T$ is the \emph{analysis
operator} $T^* : \h \to \lt$ given by
\[
T^* f = \{ \langle f , f_k \rangle \}.
\]

The \emph{frame operator} is $T T^*$. The frame operator is
particularly useful if $\ftk$ actually is a frame, but it is well
defined and bounded if we just assume $\ftk$ is a Bessel sequence.

Given a Bessel sequence $\ftk$, a Bessel sequence $\gtk$ is called a
{\it dual frame} if \bee \label{a1} f = \sum \la f,f_k \ra g_k,
\quad \forall f\in \h.\ene Condition \eqref{a1} means that analysis
using $\ftk$ followed by synthesis using $\gtk$ yields the identity
operator. We adopt the terminology of signal processing and speak
then about {\it perfect reconstruction.}

One choice of dual frame is the {\it canonical dual frame}
$\{(TT^*)^{-1}f_k\}$. When $\ftk$ is redundant, infinitely many
other dual frames exist.

Note that if $\ftk$ has the upper frame bound $B$ and $\gtk$ is a
dual frame, then an application of the Cauchy--Schwarz inequality
on \bes \lVert f\rVert ^2 = \sum \la f,f_k\ra \la g_k,f \ra \ens
shows that $\gtk$ satisfies the lower frame inequality with bound
$1/B;$ thus the sequences $\ftk$ and $\gtk$ are indeed frames. For
later use, observe that no estimate on the upper frame bound for
$\gtk$ can be deduced from the duality condition or from knowledge
about the frame bounds for $\ftk$, as Example \ref{A1} in
Appendix~\ref{framecharact_proof} shows.

\section{Various concepts of duality} \label{sec3}

In order to apply the dual frame expansion \eqref{a1} for a given
frame $\ftk$, we need to find a dual frame $\gtk$. Unfortunately,
it might be cumbersome --- or even impossible --- to calculate a
dual frame explicitly. In the literature one finds only a few
infinite dimensional non-tight frames for which a dual has been
constructed. This paucity of constructions leads us to seek frames
that are ``close to dual''. We will propose two such concepts.

Suppose that $\{ f_k \}$ and $\{ g_k \}$ are Bessel sequences in
$\h;$ we will denote their synthesis operators by $T : \lt \to \h$
and $U : \lt \to \h,$ respectively. The operators $TU^*:\h \to \h$
and $UT^*:\h \to \h$ will be called \emph{mixed frame operators.}
Note that \bes TU^*f= \sum \la f, g_k\ra f_k, \quad f\in \h,\ens and
\bes UT^*f= \sum \la f, f_k\ra g_k, \quad f\in \h.\ens Recall that
$\ftk$ and $\gtk$ are dual frames when $TU^*=I$ or $UT^*=I$.
\begin{defn} Bessel sequences $\{ f_k \}$ and $\{ g_k \}$ are said to be
\begin{itemize}
\item \textbf{approximately dual frames} if
$\lVert I-T U^*  \rVert < 1$ or $\lVert I-U T^*  \rVert < 1$,
\item \textbf{pseudo-dual frames} if $T U^*$ or $U T^*$ is a bijection on
$\h$.
\end{itemize}
\end{defn}
In each definition, the two given conditions are equivalent, by
taking adjoints. Note also that the article \cite{CHS} uses the
term approximate dual with a different meaning than we do here.

We illustrate the pseudo-dual idea by means of a frame
characterization.
\bt[Characterization of frames] \label{framecharact} Let $\{ f_k \}$
be a Bessel sequence in $\h$. Then the following statements are
equivalent.

\bei \item[(a)] $\{ f_k \}$ is a frame for $\h$. \item[(b)] The
frame operator $T T^*$ is a bijection on $\h$. \item[(c)] The
synthesis operator $T$ is surjective onto $\h$. \item[(d)] $\ftk$
has a dual frame $\gtk$. \item[(e)] $\ftk$ has a pseudo-dual frame
$\gtk$. \eni \et

The new part of the theorem is the pseudo-dual statement in part
(e). We need, however, some information from the proof of the
known parts at several instances in our paper, so we include a
full proof of Theorem \ref{framecharact} in
Appendix~\ref{framecharact_proof}.

Part (b) of the theorem shows that every frame is pseudo-dual to
itself. The proof of (a)$\Rightarrow$(b) gives even more: in fact,
denoting the frame bounds for $\ftk$ by $A$ and $B$, the result
\eqref{20a} in Appendix \ref{framecharact_proof} shows that \bee
\label{21a} \lVert I- \frac2{A+B}TT^* \rVert \le \frac{
\frac{B}{A}-1}{\frac{B}{A}+1} < 1, \ene \emph{i.e.}, every frame
$\ftk$ is approximately dual to the multiple $\{ \frac2{A+B}f_k \}$
of itself.

Next we establish relations among the three duality concepts.

\bl[Duality relations] \label{l3} Let $\ftk$ and $\gtk$ be Bessel
sequences in $\h$. \bei \item[(i)] If $\ftk$ and $\gtk$ are dual
frames, then $\ftk$ and $\gtk$ are approximately dual frames.
\item[(ii)] If $\ftk$ and $\gtk$ are approximately dual frames, then
$\ftk$ and $\gtk$ are pseudo-dual frames. \item[(iii)] If $\ftk$ and
$\gtk$ are pseudo-dual frames, then $\ftk$ and $\gtk$ are frames.
\item[(iv)] If $\ftk$ and $\gtk$ are pseudo-dual frames and $W:\h
\to \h$ is a bounded linear bijection, then $\ftk$ and $\{Wg_k\} $
are pseudo-dual frames.\eni \el \bp Statement (i) is an immediate
consequence of the definitions. So is (ii), since $\lVert I - UT^*
\rVert < 1$ implies $UT^*$ is a bijection, with inverse given by a
Neumann series. Statement (iii) follows from
Theorem~\ref{framecharact}. For the proof of (iv), observe that the
synthesis operator for $\{Wg_k\} $ is $X=WU;$ the assumptions of
$\ftk$ and $\gtk$ being pseudo-dual frames and $W$ being a bijection
imply that $XT^*=WUT^*$ is a bijection.\ep

The property of being a pair of pseudo-dual frames is significantly
weaker than being a pair of dual frames. Nevertheless, the following
proposition shows that pseudo-dual frame pairs generate dual frame
pairs in a natural fashion.
\begin{prop}[Pseudo-duals generate duals] \label{naturalduals} \
If $\{ f_k \}$ and $\{ g_k \}$ are pseudo-dual frames, then $\{
f_k \}$ and $\{ (U T^*)^{-1} g_k \}$ are dual frames. \end{prop}

\bp Assuming that $\{ f_k \}$ and $\{ g_k \}$ are pseudo-dual
frames, we know that $(U T^*)^{-1}$ exists and is bounded. Hence
$\{ (U T^*)^{-1} g_k \}$ is a Bessel sequence. If we analyze with
$\{ f_k \}$ and synthesize with $\{ (U T^*)^{-1} g_k \}$ then we
obtain the identity, since
\begin{align*}
\sum \langle f , f_k \rangle (U T^*)^{-1} g_k & = (U T^*)^{-1}
\sum \langle f , f_k \rangle g_k  \\
& = (U T^*)^{-1} U T^* f \\
& = f .
\end{align*}
Thus $\{ f_k \}$ and $\{ (U T^*)^{-1} g_k \}$ are dual frames. \ep

The relation of being a pair of pseudo-dual frames is symmetric: if
$\{ f_k \}$ and $\{ g_k \}$ are pseudo-dual frames, then so are $\{
g_k \}$ and $\{ f_k \}$. The pseudo-dual relation is also reflexive
on the set of frames since every frame is pseudo-dual to itself, as
remarked after Theorem~\ref{framecharact}. Transitivity of the
pseudo-dual relation fails (in Hilbert spaces of dimension at least
2) by Example~\ref{r1} below.

\section{Approximately dual frames} \label{sec4}

In this section we focus on approximately dual frames $\ftk$ and
$\gtk.$ As before, we denote the associated synthesis operators by
$T$ and $U$, respectively. Our main goal is to demonstrate how one
can use this concept to obtain what in engineering terms would be
called ``almost perfect reconstruction''.

Recall that if $\ftk$ and $\gtk$ are approximately dual frames
then \bes \lVert  f- \sum \la f,f_k \ra g_k  \rVert = \lVert
(I-UT^*)f \rVert \le \lVert I-UT^* \rVert \, \lVert f \rVert.\ens
If $\lVert I-UT^* \rVert <<1$ then we see that $\sum \la f,f_k \ra
g_k$ ``almost reconstructs'' $f$, which motivates our terminology.

We first strengthen Proposition \ref{naturalduals}, for
approximately dual frames and the associated ``natural'' dual frame.

\begin{prop} \label{na} Assume that $\{ f_k \}$ and $\{ g_k \}$ are
approximately dual frames. Then the following hold: \bei \item[(i)]
The dual frame $\{(U T^*)^{-1} g_k\}$ of $\ftk$ can be written \bee
\label{a5} (U T^*)^{-1} g_k  = g_k + \sum_{n=1}^\infty (I - UT^*)^n
g_k . \ene
\item[(ii)] Given $N\in \mn$,
consider the corresponding partial sum, \bee \label{a6}
\ga_k^{(N)}= g_k + \sum_{n=1}^N (I - UT^*)^n g_k.\ene Then
$\{\ga_k^{(N)}\}$ is an approximate dual of $\ftk$. Denoting its
associated synthesis operator by $Z$, we have \bee \label{12}
\lVert I-ZT^* \rVert \le  \lVert I-UT^* \rVert^{N+1}.\ene \eni
\end{prop}
\bp Assuming that $\ftk$ and $\gtk$ are approximately dual frames,
the inverse of $UT^*$ can be written via a Neumann series as \bee
\label{7a} (UT^*)^{-1} = \big( I - (I - UT^*) \big)^{-1} =
\sum_{n=0}^\infty (I - UT^*)^n . \ene The result in (i) now
follows by applying this expansion to $g_k$ and recalling
Proposition \ref{naturalduals}.

For (ii), note that $\{\ga_k^{(N)}\}$ is a Bessel sequence since
it is obtained from the Bessel sequence $\gtk$ by a bounded
transformation. And \bes ZT^*f & =
& \sum \la f,f_k\ra \big (I+ \sum_{n=1}^N (I - UT^*)^n \big)g_k \\
& = & \sum_{n=0}^N (I - UT^*)^n UT^*f \\
& = & \sum_{n=0}^N (I - UT^*)^n \big( I - (I-UT^*) \big)f \\
& = & f - (I-UT^*)^{N+1} f\ens by telescoping. Thus \bee \notag
\lVert I-ZT^* \rVert & = & \lVert (I-UT^*)^{N+1} \rVert \\
\notag & \leq & \lVert I-UT^* \rVert^{N+1} < 1.\ene
 \ep

\begin{rem} \rm In view of formula \eqref{a5} in Proposition \ref{na},
an approximately dual frame $\{ g_k \}$ associated with a frame
$\ftk$ can be regarded as a zero-th order approximation to the
(exact) dual frame $\{(U T^*)^{-1} g_k\}$.  In case $\lVert I-UT^*
\rVert$ is small, reconstruction using the approximate dual $\gtk$
is close to perfect reconstruction. The result in (ii) yields a
family of approximately dual frames that interpolates between the
approximate dual $\gtk$ and the dual frame $\{(U T^*)^{-1} g_k\}$;
the estimate \eqref{12} shows that by choosing $N$ sufficiently
large, we can obtain a reconstruction that is arbitrarily close to
perfect. The drawback of the result with respect to potential
applications is the complicated structure of the operator in
\eqref{a6} defining the sequence $\ga_k^{(N)}$.
\end{rem}

We have now presented the basic facts for approximately dual frames.
In the following, we will present concrete settings where they yield
interesting new insights.

To motivate the results, we note (again) that it can be a
nontrivial task to find the canonical dual frame (or any other
dual) associated with a general frame $\ftk.$ We seek to connect
this fact with perturbation theory by asking the following
question: if we can find a frame $\htk$ that is close to $\ftk$
and for which it is possible to find a dual frame $\{g_k\}$, does
it follow that $\{g_k\}$ is an approximate dual of $\ftk$? We will
present some sufficient conditions for an affirmative answer.

First we state a general result, valid for dual frames $\gtk$ with
sufficiently small Bessel bound; later we state a more explicit
consequence for the case where $\gtk$ is the canonical dual frame,
see Theorem \ref{c1}.  Technically, we need not assume that $\ftk$
is a frame in the general result: the frame property follows as a
conclusion. On the other hand, as explained above, the main use of
the result is in a setting where $\ftk$ is known to be a frame in
advance.

\bt[Dual of perturbed sequence] \label{t1} Assume that $\ftk$ is a
sequence in $\h$ and that $\htk$ is a frame for which \[ \sum | \la
f,f_k-h_k\ra|^2 \le R \, \lVert f \rVert^2, \quad \afh,\] for some
$R>0$. Consider a dual frame $\gtk$ of $\htk$ with synthesis
operator $U$, and assume $\gtk$ has upper frame bound $C$.

If $CR<1$ then $\ftk$ and $\gtk$ are approximately dual frames, with
\[ \lVert I-UT^* \rVert \le \lVert U \rVert \, \sqrt{R} \le \sqrt{CR} <
1.\] \et

\bp With our usual notation, we have $UV^*=I$ since $\gtk$ and
$\htk$ are dual frames. Hence \bes \lVert I- UT^* \rVert = \lVert
U(V^*-T^*) \rVert \le \lVert U \rVert \, \lVert V^* - T^* \rVert \le
\sqrt{CR} <1.\ens \ep

It is crucial in Theorem \ref{t1} that the dual frame $\gtk$ has
upper frame bound less than $1/R$. Otherwise $\gtk$ need not be
pseudo-dual to $\ftk$, let alone approximately dual, as the next
example shows.

\bex \label{r1} \rm Consider $\h = \mc^2$ with the standard basis
$\{ e_1 , e_2 \}$. Let $\epsilon > 0$ and consider the frames
\[ \ftk = \{0, e_1, e_2\}, \quad \htk = \{\epsilon e_1, e_1, e_2\}, \quad
\gtk= \{ \epsilon^{-1} e_1, 0, e_2\} . \] Write $T,V,U$ for the
associated synthesis operators.  Note that \bes T U^* f = \langle f ,
e_2 \rangle e_2;\ens this operator is neither injective nor
surjective, and so $\{ f_k \}$ and $\{ g_k \}$ are not pseudo-dual
frames. Clearly, $\htk$ is a frame for $\mc^2,$ regardless of the
choice of $\epsilon.$ Now, because \bes \sum_{k=1}^3 | \la f,
f_k-h_k\ra|^2 = |\la f, \epsilon e_1\ra|^2 \le \epsilon^2 \lVert f
\rVert^2 , \quad \forall f\in \mc^2 ,\ens the condition in Theorem
\ref{t1} is satisfied with $R=\epsilon^2.$ Thus we see that no
matter how close $\htk$ gets to $\ftk$ (meaning, no matter how
small $\epsilon$ is), the frame $\gtk$ is not pseudo-dual to
$\ftk$. Note that for $\epsilon<1,$ the frame $\gtk$ has the upper frame bound
$C=\epsilon^{-2}=1/R;$ thus,
Theorem \ref{t1} is not contradicted.

In this example it is immediate to see that $T V^* = I$ and $V U^*
= I;$ thus $\{ f_k \}$ and $\{ h_k \}$ are dual frames, and so are
$\{ h_k \}$ and $\{ g_k \}$. In particular, the example shows that
the pseudo-dual property is not transitive. The same argument
applies in any separable Hilbert space of dimension at least 2.
 \ep \enx

For practical applications of Theorem~\ref{t1}, it is problematic
that the Bessel bound $C$ on $\gtk$ might increase when one
reduces $R$ by taking $\htk$ very close to $\ftk:$ in this way, it
might not be possible to get $CR<1$ just by taking $\htk$
sufficiently close to $\ftk$. In the case when $\gtk$ is the
\emph{canonical} dual of $\htk$ the problem does not arise: as the
next result shows, if $\htk$ is sufficiently close to $\ftk$, then
the canonical dual $\gtk$ of $\htk$ is approximately dual to
$\ftk$. For the statement of this result we need to assume that
$\ftk$ is a frame, since its lower frame bound controls the
perturbation hypothesis.

\bt[Canonical dual of perturbed frame] \label{c1} Let $\ftk$ be a
frame for $\h$ with frame bounds $A,B.$ Let $\htk$ be a sequence
in $\h$ for which \bes \sum | \la f,f_k-h_k\ra|^2 \le R \, \lVert
f \rVert^2, \quad \afh,\ens for some $R < A/4.$ Denote the
synthesis operator for $\htk$ by $V.$

Then $\htk$ is a frame. Its canonical dual frame
$\gtk=\{(VV^*)^{-1}h_k\}$ is an approximate dual of $\ftk$ with
\bes \lVert I-UT^* \rVert  \le \frac{1}{\sqrt{A/R} - 1} < 1 , \ens
where $U$ denotes the synthesis operator for $\{g_k\}$. \et

\bp The sequence $\htk$ is a frame with frame bounds
$(\sqrt{A}-\sqrt{R})^2$ and $(\sqrt{B}+\sqrt{R})^2$, by
\cite[Corollary 5.6.3]{CBN}; the proof is short, and so we give it
here. We have $\lVert T^*-V^* \rVert^2 \le R$ by assumption, and
so
\[ | \lVert T^* f \rVert_\lt - \lVert V^* f \rVert_\lt | \leq
\sqrt{R} \lVert f \rVert \] for all $f \in \h$, from which the
desired frame bounds follow easily. This part of the proof uses
only $R<A$.

The canonical dual frame of $\htk$ is $\gtk = \{(VV^*)^{-1}h_k\},$
with frame bounds \bes \frac1{(\sqrt{B} + \sqrt{R})^2} \quad
\text{and} \quad \frac1{(\sqrt{A} - \sqrt{R})^2} \ens by
\cite[Lemma 5.1.6]{CBN}. In terms of the synthesis operator $U$
for $\gtk$, the upper bound says \[\lVert U \rVert \leq
\frac1{\sqrt{A}-\sqrt{R}} .\]

Theorem \ref{t1} therefore implies \bee \notag
\lVert  I-UT^* \rVert & \le & \lVert U \rVert \, \sqrt{R} \\
& \le & \frac{\sqrt{R}}{\sqrt{A} - \sqrt{R}} \notag \\
\notag & = & \frac{1}{\sqrt{A/R} - 1}.\ene To complete the proof,
just notice this last expression is smaller than $1$ if and only
if $R< A/4.$ \ep

We end this section with remarks concerning the appearance of
approximate duals in the literature. Approximately dual frames have
been employed in the $L^2$ theory of wavelets; for example,
Holschneider developed a sufficient condition for a pair of wavelet
systems to be approximately dual \cite[Section 2.13]{H}, and Gilbert
\emph{et al.}\ \cite{GHHLWW} obtained approximate duals for highly
oversampled wavelet systems by perturbing the continuous parameter
system. Approximate duals have proved central to the recent solution
of the Mexican hat wavelet spanning problem \cite{BL}. There the
synthesizing wavelet frame is given, being generated by the Mexican
hat function, and the task is to construct an approximately dual
analyzing frame. Interestingly, this approximate dual construction
holds in all $L^q$ spaces (in the frequency domain), which suggests
that parts of the theory in this paper might extend usefully to
Banach spaces.

Approximately dual frames in Gabor theory arose in the work of
Feichtinger and Kaiblinger \cite[Sections 3,4]{FK}. They studied the
stability of Gabor frames with respect to perturbation of the
generators and the time-frequency lattice; among many other results,
they showed that a sufficiently small perturbation (measured in a
specific norm) of a dual frame associated with a Gabor frame leads
to an approximately dual frame. In the next sections, we construct
\emph{explicit} examples of approximately dual Gabor frames.

\section{Gabor frames and approximate duals} \label{gab}
The general frame estimates in Section \ref{sec4} are based on
operator inequalities such as the product rule $\lVert TU \rVert \le
\lVert T \rVert \lVert U \rVert$. We will apply these general
estimates in the next section to obtain explicit approximate duals
for Gabor frames generated by a Gaussian. First, though, we will
sharpen the general frame estimates in the concrete, rich setting of
Gabor theory. These sharper estimates also will be applied in the
next section.

A {\it Gabor frame} is a frame for $\ltr$ of the form $$\{e^{2\pi
imb}g(x-na)\}_{m,n\in \mz}$$ for suitably chosen parameters $a,b>0$
and a fixed function $g\in \ltr,$ called the {\it window function.}
The number $a$ is called the {\it translation parameter} and $b$ is
the {\it modulation parameter.} Introducing the operators \bes (T_a
g)(x) = g(x-a), \ (E_b g)(x)=e^{2\pi ibx}g(x), \quad \text{for\
}a,b,x\in \mr,\ens the Gabor system can be written in the short form
$\mts.$ For more information on Gabor analysis and its role in
time--frequency analysis we refer to the book by Gr\"ochenig
\cite{G2}.

We state now the duality conditions for a pair of Gabor systems, due
to Ron and Shen \cite{RoSh5}. We will apply the version presented by
Janssen \cite{J}:

\bl \label{l5} Two Bessel sequences
$\{E_{mb}T_{n}\varphi\}_{m,n\in \mz}$ and
$\{E_{mb}T_{n}g\}_{m,n\in \mz}$ form dual frames for $\ltr$ if and
only if the equations \bee \label{gfs}
\sum_{k\in \mz} \overline{\varphi (x-ak)}g(x-ak)- b & = & 0, \\
\label{gfs2} \sum_{k\in \mz} \overline{\varphi(x-n/b-ak)}g(x-ak) & =
& 0,  \quad n\in \mz \setminus \{0\} , \ene hold a.e. \el

Recall that the Wiener space $W$ consists of all bounded measurable
functions $g:\mr \to \mc$ for which \bes \sum_{n\in \mz} \lVert g
\chi_{[n,n+1[} \rVert_\infty < \infty.\ens It is well known that if
$g\in W$ then $\mts$ is a Bessel sequence for each choice of
$a,b>0.$

\bt \label{t4} Given two functions $\varphi, g \in W$ and two
parameters $a,b>0$, let $T$ denote the synthesis operator associated
with the Gabor system $\{E_{mb}T_{na}\varphi\}_{m,n\in \mz},$ and
$U$ the synthesis operator associated with
$\{E_{mb}T_{na}g\}_{m,n\in \mz}.$ Then \bes \lVert I-UT^* \rVert \le
\frac1{b} \left[ \Big\lVert b- \sum_{k\in \mz}
\overline{T_{ak}g}T_{ak}\varphi \Big\rVert_\infty+\sum_{n\neq 0}
\Big\lVert \sum_{k\in \mz} \overline{T_{n/b}T_{ak}g}T_{ak}\varphi
\Big\rVert_\infty \right].\ens \et

\bp The starting point is the Walnut representation of the mixed
frame operator associated with $\{E_{mb}T_{na}\varphi\}_{m,n\in
\mz}$ and $\mts.$   According to Theorem 6.3.2 in \cite{G2},

\bes UT^*f(\cdot)= \frac1{b} \sum_{n\in \mz} \left(\sum_{k\in \mz}
\overline{T_{ak}\varphi( \cdot-n/b)}T_{ak}g(\cdot) \right)
T_{n/b}f(\cdot), \quad f\in \ltr. \ens Thus, \bes & \ & \lVert  f-UT^*f \rVert \\
\notag & \le & \Big\lVert  \Big( 1-\frac1{b}\sum_{k\in \mz}
\overline{T_{ak}\varphi(\cdot)}T_{ak}g(\cdot) \Big) f \Big\rVert
\\ \notag & \ &
+ \frac1{b} \Big\lVert  \sum_{n\neq 0} \Big(\sum_{k\in \mz}
\overline{T_{ak}\varphi(\cdot-n/b)}T_{ak}g(\cdot) \Big) T_{n/b}f
\Big\rVert
\\ & \notag \le & \frac1{b} \Big\lVert b- \sum_{k\in \mz}
\overline{T_{ak}\varphi}T_{ak}g \Big\rVert_\infty \, \lVert f\rVert
+ \frac1{b} \sum_{n\neq 0} \Big\lVert \sum_{k\in \mz}
\overline{T_{n/b}T_{ak}\varphi}T_{ak}g \Big\rVert_\infty \, \lVert
f\rVert , \label{30a} \ens which concludes the proof.\ep Observe
that the terms appearing in the estimate in Theorem \ref{t4} measure
the deviation from equality in the duality relations in Lemma
\ref{l5}, with respect to the $\lVert \cdot\rVert _\infty$-norm. In
particular, Theorem \ref{t4} proves the ``sufficient'' direction of
Lemma \ref{l5}, under the stated hypotheses.

\section{Applications to Gabor frames generated by the Gaussian}
\label{sec16}

The Gaussian $g(x)=e^{-x^2}$ is well known to generate a Gabor frame
whenever $ab<1.$ For the case where the parameters $a$ and $b$ are
equal and small, Daubechies has demonstrated that the Gabor frame
generated by the Gaussian is almost tight, see \cite{Da2}, p.84--86
and \cite{Da1}, p. 980--982; in particular, the formula \eqref{21a}
yields almost perfect reconstruction.

Regardless of the choice of $a$ and $b$, no convenient explicitly
given expression is known for any of the dual frames associated with
the Gaussian, though. We will construct explicit, approximately dual
frames associated with the Gaussian, for certain choices of $a$ and
$b$. These approximately dual frames provide \emph{almost} perfect
reconstruction; Example \ref{22a} is particularly interesting
because it deals with a frame that is far from being tight,
\emph{i.e.}, no easy way of obtaining an approximately dual frame is
available.

\bex \label{e1} Consider the (scaled) Gaussian \bee \label{ds}
\varphi(x)= \frac{151}{315}e^{-(x/1.18)^2}.\ene The scaling is
introduced for convenience; similar constructions can be performed
for other Gaussians as well. It is well known that for this
Gaussian, the Gabor system $\{E_{mb}T_n\varphi\}_{m,n\in \mz}$ forms
a frame for $\ltr$ for any sufficiently small value of the
modulation parameter $b>0$; we fix $b=0.06$ in this example. Denote
the synthesis operator for that frame by $T$.

\begin{figure}
\centerline{
\includegraphics[width=5in,height=4in]{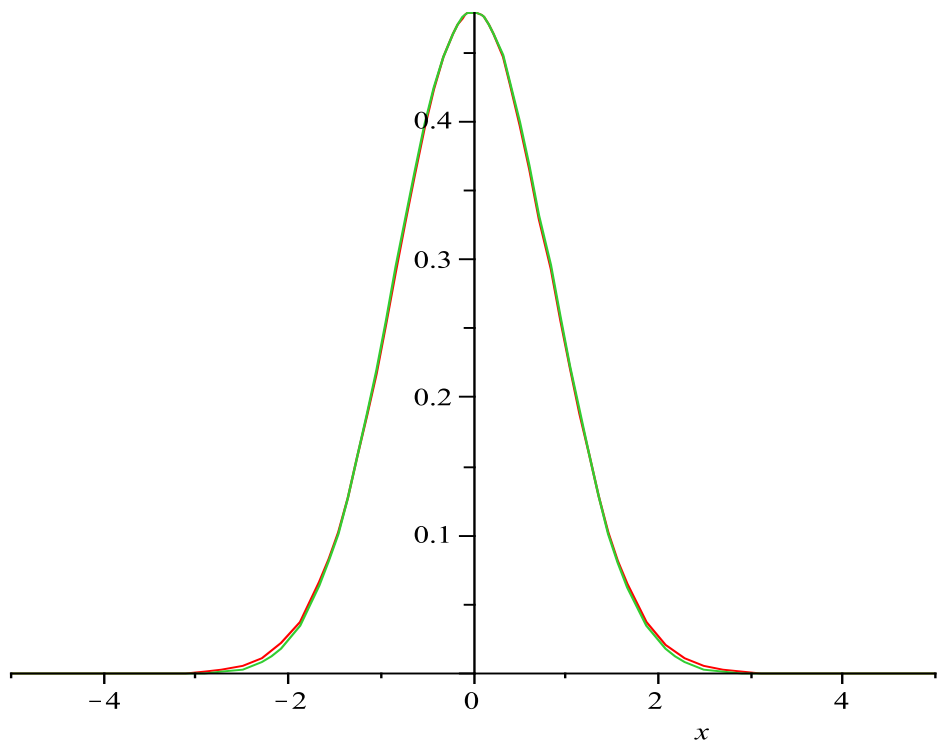}}
\caption{The B-spline $B_8$ and the scaled Gaussian $\varphi$ in
\eqref{ds}.} \label{ff1}
\end{figure}

\begin{figure}
\centerline{
\includegraphics[width=2.7in,height=2.5in]{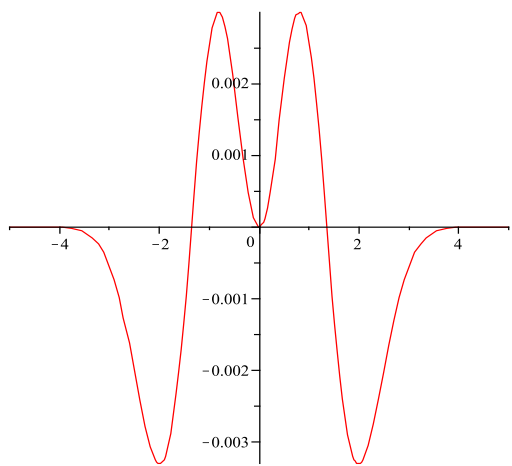}}
\caption{The function $B_8-\varphi.$} \label{ff2}
\end{figure}

\begin{figure}
\centerline{
\includegraphics[width=2.7in,height=2.5in]{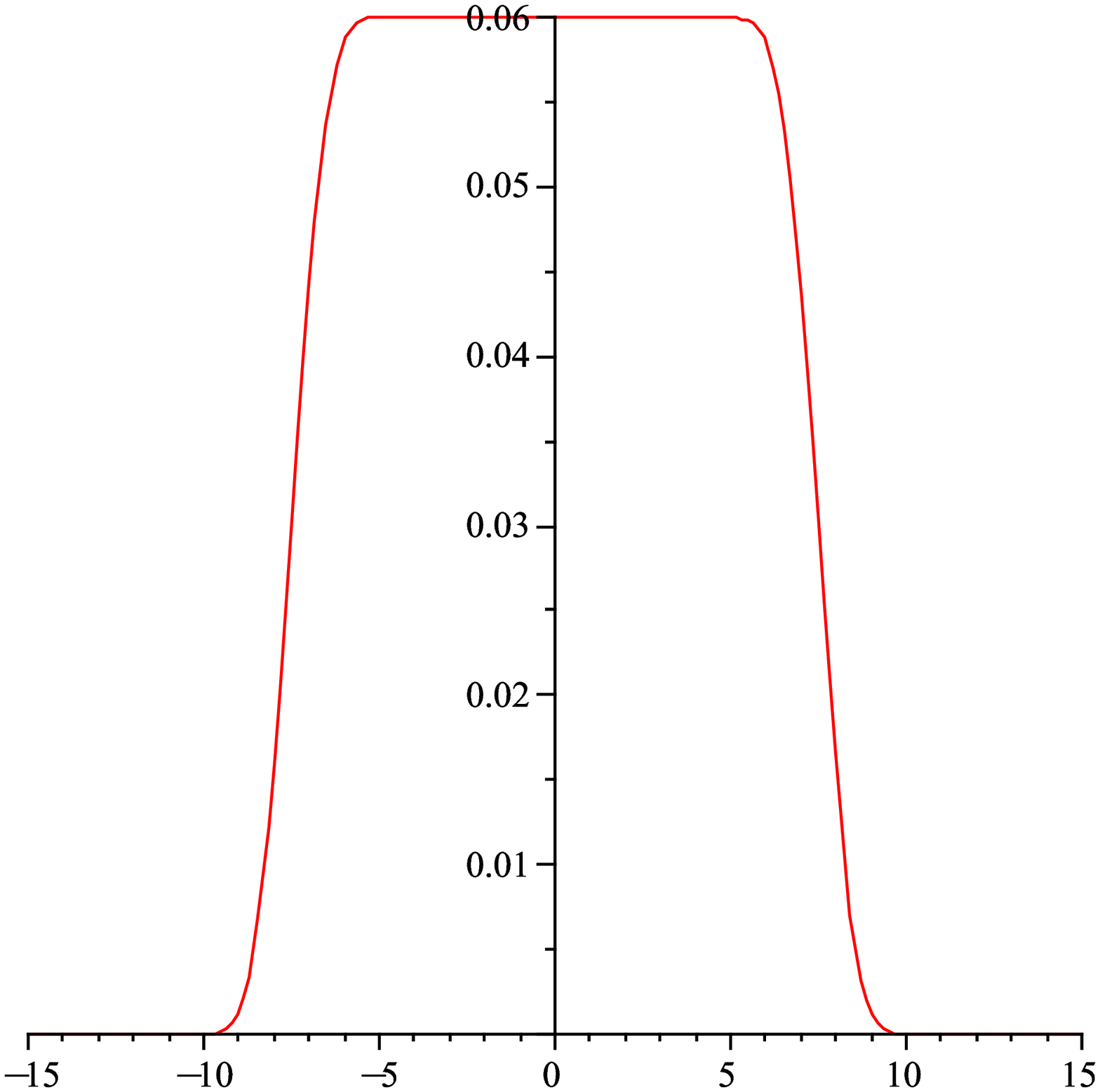}}
\caption{The non-canonical dual window $g$ of $B_8$, see
\eqref{2a}.} \label{ff3}
\end{figure}

We will use the results derived in this paper to find an
approximately dual Gabor frame. Let $h=B_8$ denote the
eighth-order B-spline, centered at the origin. This B-spline
approximates the Gaussian $\varphi$ very well: see
Figures~\ref{ff1} and \ref{ff2}. The Gabor system $\{E_{b
m}T_nh\}_{m,n\in \mz}$ is a frame for $\ltr$, and by Corollary 3.2
in \cite{CK} it has  the (non-canonical) dual frame $\{E_{b
m}T_ng\}_{m,n\in \mz},$ where  \bee \label{2a} g(x)=b
\sum_{n=-7}^7 B_8(x+n)\ene is a linear combination of translated
B-splines. This $g$ is supported on $[-11,11]$ and is constantly
equal to $b$ on $\text{supp} \, B_8=[-4,4];$ see Figure~\ref{ff3}.

We first apply Theorem \ref{t1}. A numerical calculation based on
Theorem 9.1.5 in \cite{CBN} shows that the functions $\{E_{b
m}T_n(\varphi-h)\}_{m,n\in \mz}$ form a Bessel sequence with
Bessel bound $R=0.0006.$

Furthermore, a second application of Theorem 9.1.5 in \cite{CBN}
reveals that $\{E_{b m}T_ng\}_{m,n\in \mz}$ has Bessel bound $C=1.$
Thus the hypotheses of Theorem \ref{t1} are satisfied. Denoting the
synthesis operator for $\{E_{b m}T_n g\}_{m,n\in \mz}$ by $U$, we
deduce that \bes \lVert I-UT^* \rVert \le \sqrt{CR}\le 0.025.\ens On
the other hand, an application of Theorem \ref{t4} immediately
yields the better estimate \bes \lVert I-UT^* \rVert \le 0.0027.\ens

In terms of the Gabor systems involved, this conclusion says that
\bes \lVert f- \sum_{m,n\in \mz} \la f, E_{b m}T_n \varphi \ra E_{b
m}T_ng \rVert < 0.0027 \, \lVert f \rVert_2, \quad \forall f \in
\ltr .\ens That is, analysis with the scaled Gaussian and synthesis
with the function $g$ yields almost perfect reconstruction.

We note that the frame considered in this example is almost tight;
thus, the approximate dual constructed here does not perform better
than using a scaled version of the frame itself as approximate dual.
For a frame that is far from being tight this aspect changes
drastically in favor of the approximate duals considered in this
paper -- see Example \ref{22a}. Lastly, we remark that other
approximations to the Gaussian are possible too; instead of using
B-splines and their duals from \cite{CK} we could have used certain
other splines and their duals from \cite{L}.

\ep \enx

In the next example we apply Theorem \ref{c1} and Theorem
\ref{t4}. Let us again consider  B-spline approximation to a
Gaussian. The advantage of using the canonical dual and
Theorem~\ref{c1}, compared with our use of a non-canonical dual
frame in Example \ref{e1}, is that larger values of the modulation
parameter $b$ can be handled. Approximating using $B_8$, our
approach in Example \ref{e1} is restricted to $b\le 1/15$ (the
restriction comes from the underlying results in \cite{CK}); as
illustration of the larger range for $b$ obtained via Theorem
\ref{c1},  we take $b=0.1$ in the next example.

\bex \label{22a} Let \bes \varphi(x)= e^{-4x^2}.\ens An
application of Theorem 5.1.5 in \cite{CBN} shows that the functions
$\{E_{0.1m}T_n\varphi\}_{m,n\in \mz}$ form a frame for $\ltr$ with
frame bounds $A=2.6, \ B=10.1$. The function
\bes h(x)= \frac{315}{151}B_8(2.36x)\ens
yields a close approximation of $\varphi,$ see Figure 4.
A numerical calculation
based on Theorem 5.1.5 in \cite{CBN}  shows that the functions
$\{E_{0.1 m}T_n(\varphi-h)\}_{m,n\in \mz}$ form a Bessel sequence with
Bessel bound $R=6.5\cdot 10^{-4} <A/4.$

The function $h$ has support on $[-4/2.36,4/2.36],$ an
interval of length $8/2.36$. Since the
modulation parameter $b=0.1$ is smaller than
$(8/2.36)^{-1}$  and the function
\bes H(x)= \sum_{k\in \mz} |h(x+k)|^2\ens is bounded above and below away
from zero,
it follows from Corollary 9.1.7 in \cite{CBN} that
the frame
$\{E_{0.1m}T_nh\}_{m,n\in \mz}$ has the canonical dual frame
$\{E_{0.1m}T_ng\}_{m,n\in \mz},$  where \bee \notag g(x) & = &
\frac{0.1}{\sum_{n\in \mz} |h(x+n)|^2} \, h(x) \\ \label{19a} & = &
\frac{15.1}{315} \, \frac1{\sum_{n\in \mz} |B_8(2.36(x+n))|^2} \,
B_8(2.36x).
\ene See Figure 5.

\begin{figure}
\centerline{
\includegraphics[width=2.7in,height=2.5in]{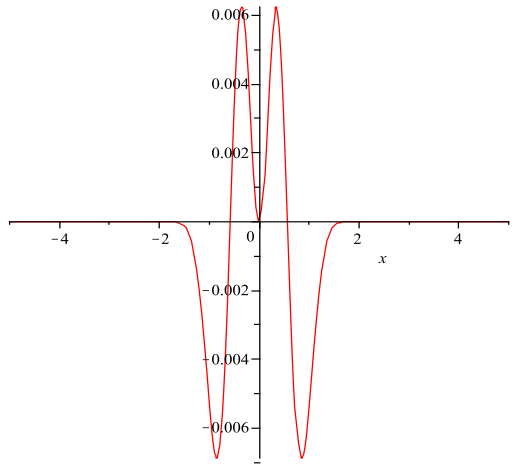}}
\caption{The function $x\mapsto e^{-4x^2}-\frac{315}{151}B_8(2.36x). $}
\end{figure}

\begin{figure}
\centerline{
\includegraphics[width=2.7in,height=2.5in]{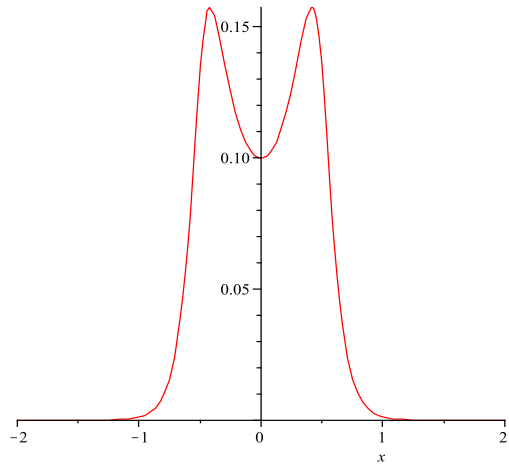}}
\caption{The approximative dual window $g$ in \eqref{19a}.}
\end{figure}

The frame $\{E_{0.1m}T_n g\}_{m,n\in \mz}$ is approximately dual to
$\{E_{0.1m}T_n\varphi\}_{m,n\in \mz}$, by Theorem \ref{c1}. Denoting
their synthesis operators by $U$ and $T$ respectively, the
approximation rate in Theorem \ref{c1} is measured by \bee
\label{31a} \lVert I-UT^* \rVert \le \frac{1}{\sqrt{A/R} - 1} \le
0.016.\ene On the other hand, Theorem \ref{t4} yields the somewhat
better estimate \bee \label{31af} \lVert I-UT^* \rVert \le
\frac{1}{\sqrt{A/R} - 1} \le 0.009.\ene Thus the approximate dual
frame almost yields perfect reconstruction.

Note that the frame $\{E_{0.1m}T_n\varphi\}_{m,n\in \mz}$ is far
from being tight. Thus, the sequence $\{E_{0.1m}T_n
\left(\frac2{A+B}\, \varphi\right)\}_{m,n\in \mz}$ is a poor
approximate dual: the estimate in \eqref{21a} yields \bes \lVert
I- \frac2{A+B}S \rVert \le \frac{\frac{B}{A}-1}{ \frac{B}{A}+1} =
0.59,\ens which is far worse than the result in \eqref{31af}. \ep
\enx

\begin{rem}\rm \

\noindent 1. In case a closer approximation than in \eqref{31af} is
required, an application of the iterative procedure in Proposition
\ref{na} (ii) can bring us as close to perfect reconstruction as
desired. Letting $\ftk$ be a Gabor frame $\{E_{mb}T_{na}f\}_{m,n\in
\mz}$ and $\gtk$ be an approximately dual Gabor frame $\mts$, each
approximate dual $\{\ga_k^{(N)}\}$ constructed in Proposition
\ref{na} will again have Gabor structure. For $N=1,$ an easy
calculation shows that $\{\ga_k^{(N)}\}$ will be the Gabor system
$\{E_{mb}T_{na}\gamma\}_{m,n\in \mz},$ where \bes \gamma =
g+(I-UT^*)g & = & 2g- UT^*g \\ & = & 2g - \sum_{m^\prime,n^\prime
\in \mz}
\la g, E_{m^\prime b}T_{n^\prime a}f \ra E_{m^\prime b}T_{n^\prime a} g \\
& = & ( 2- \la g, f\ra) \, g - \sum_{(m^\prime,n^\prime)\neq (0,0)}
\la g, E_{m^\prime b}T_{n^\prime a}f \ra E_{m^\prime b}T_{n^\prime
a} g.\ens For the particular case considered in Example \ref{22a},
and denoting the synthesis operator for
$\{E_{mb}T_{na}\gamma\}_{m,n\in \mz}$ by $Z,$ Proposition \ref{na}
shows that the estimate in \eqref{31a} will be replaced by \bes
\lVert I-Z^*T\rVert  \le 0.009^2= 8.1 \times 10^{-5}.\ens

\noindent 2. Another way of calculating approximate duals for Gabor
frames based on $\varphi(x)=e^{-4x^2}$ would be to follow the
procedure in Example \ref{22a} with a function $h$ of the type \bes
h(x)= \varphi(x) \chi_I(x)\ens for a sufficiently large interval
$I$, rather than letting $h$ be a B-spline. However, for such a
function $h$ the canonical dual window associated with
$\{E_{0.1m}T_n\varphi\}_{m,n\in \mz}$ will not be continuous; thus,
it will have bad time-frequency properties.
\end{rem}

\section*{Acknowledgments} \label{sec6}

Laugesen thanks the Department of Mathematics at the Technical
University of Denmark for travel support and a warm welcome.

\appendix

\section{Proofs and additional examples} \label{framecharact_proof}
In this Appendix we collect some additional information related to
classical frame theory.

First, we present the example announced in Section~\ref{sec2},
showing  that an estimate on the upper frame bound for a frame
$\gtk$ neither can be deduced from the fact that it is dual to a
frame $\gtk$, nor from knowledge of the frame bounds for $\ftk$.

\bex \label{A1} Consider the Hilbert space $\h = \mc^2$ with the
standard orthonormal basis $\{ e_1 , e_2 \}$. Let \bes \{ f_1 ,
f_2 , f_3 \} = \{ 0 , e_1 , e_2 \}.\ens  Then $\ftk$ is a frame
with bounds $A=B=1$. For any $C\in \mr,$ the sequence \bes \{ g_1
, g_2 , g_3 \} = \{ Ce_1 , e_1 , e_2 \}\ens is a dual frame of
$\ftk.$ The (optimal) upper frame bound for $\gtk$ is $C^2+1,$
which can be arbitrarily large.

Exactly the same considerations apply in an arbitrary separable
Hilbert space of dimension at least two. \ep \enx

We now state the announced proof of Theorem \ref{framecharact}.
The equivalences between the first four statements are well known;
we include the proof in order to keep the paper self-contained,
and because we refer to some of the steps elsewhere in the paper.

\vspace{.1in}\noindent {\bf Proof of Theorem \ref{framecharact}:}
First recall that $T$ and $U$ are bounded operators from $\lt$ to
$\h$, whenever $\{ f_k \}$ and $\{ g_k \}$ are Bessel sequences.

(a) implies (b): First, write $S = T T^*$ for the frame operator,
and note it is selfadjoint. By part (a) we have $A \lVert f \rVert^2
\leq \langle Sf , f \rangle \leq B \lVert f \rVert^2$ for all $f \in
\h$. That is, $AI \leq S \leq BI$; following the argument by
Gr\"ochenig \cite{G2} p.91, we infer that $ - \frac{B-A}{B+A} I \leq I -
\frac{2}{A+B} \, S \leq \frac{B-A}{B+A}I$. Thus, \bee \label{20a}
\lVert I - \frac{2}{A+B}S \rVert = \sup_{\lVert f \rVert
= 1} | \langle (I-\frac{2}{A+B} S)f , f \rangle | & \leq & \frac{B-A}{B+A}  \\
\notag & < & 1.\ene Therefore $\frac2{A+B}S$ is invertible, so that
$S$ is a bijection as desired.

(b) implies (c) is immediate.

(c) implies (d): consider the pseudo-inverse $T^\dagger : \h \to
\lt$. Applying the Riesz representation theorem to the functional
$(T^\dagger f)_k =$ ($k$-th component of $T^\dagger f$), we obtain
some element $g_k \in \h$ satisfying $(T^\dagger f)_k = \langle f ,
g_k \rangle$. The sequence $\{ g_k \}$ is Bessel, since
\[
\sum | \langle f , g_k \rangle |^2 = \sum | (T^\dagger f)_k |^2 =
\lVert T^\dagger f \rVert_{\lt}^2 \leq \lVert T^\dagger
\rVert_{\lt}^2 \lVert f \rVert^2 .
\]
The analysis operator $U^*$ associated to the sequence $\{ g_k \}$
equals exactly the pseudo-inverse $T^\dagger$, by construction, and
so $T U^* = T T^\dagger$ equals the identity, by the fundamental
property of the pseudo-inverse. Thus $\gtk$ is a dual frame for
$\ftk$, and so (d) holds.

(d) implies (e) trivially.

(e) implies (a) as follows: suppose $\gtk$ is a pseudo-dual frame
for $\ftk$, so that $TU^*$ is a bijection on $\h$. Then for all $f
\in \h$,
\begin{align*}
\lVert f \rVert
& = \lVert (T U^*)^{-1} T U^* f \rVert \\
& \leq \lVert (T U^*)^{-1} T \rVert \lVert U^* f \rVert_{\lt} ,
\end{align*}
and of course $\lVert U^* f \rVert_{\lt} \leq \lVert U^* \rVert
\lVert f \rVert$. Hence $\ftk$ is a frame with bounds $A=1/\lVert (T
U^*)^{-1} T \rVert^2$ and $B=\lVert U^* \rVert^2$. \ep

\vspace{.1in} \noindent{\it   Ole Christensen, Department of
Mathematics, Technical University of
\\ Denmark,
Building 303, 2800 Lyngby, Denmark \quad
Ole.Christensen\@@mat.dtu.dk

\vspace{.1in} \noindent Richard S. Laugesen, Department of
Mathematics, University of Illinois at Urbana--Champaign, Urbana, IL
61801, U.S.A. \quad Laugesen\@@illinois.edu }


\begin{thebibliography}{10}
\baselineskip 14pt

\bibitem{BL} H.-Q. Bui and R. S. Laugesen: {\it Frequency-scale frames and
the solution of the Mexican hat problem.} Preprint,  2008. \\
{\small \url{www.math.uiuc.edu/~laugesen/}}

\bibitem{CBN} O. Christensen: {\it Frames and bases. An
introductory course.} Birkh\"{a}user, Boston, 2008.

\bibitem{CK} O. Christensen and R. Y. Kim: {\it On dual Gabor
frame pairs generated by polynomials.}  Accepted for publication
in J. Fourier Anal. Appl., 2008.

\bibitem{CHS} Chui, C., He, W., and St\"ockler, J.: {\it
Nonstationary tight wavelet frames, I: Bounded intervals.} Appl.
Comp. Harm. Anal. {\bf 17} (2004), 141--197.

\bibitem{Da2} I. Daubechies: {\it Ten lectures on wavelets.}
SIAM, Philadelphia, 1992.

\bibitem{Da1}  I. Daubechies: {\it The wavelet transformation,
time-frequency localization and signal analysis.} IEEE Trans.
Inform. Theory {\bf 36} (1990), 961--1005.

\bibitem{FK} H. G. Feichtinger and N. Kaiblinger.
{\it Varying the time-frequency lattice of Gabor frames.} Trans.
Amer. Math. Soc. {\bf 356}, 2001--2023, 2004.

\bibitem{GHHLWW}
J. E. Gilbert, Y. S. Han, J. A. Hogan, J. D. Lakey, D. Weiland and
G. Weiss. {\it Smooth molecular decompositions of functions and
singular integral operators.}  Mem. Amer. Math. Soc. {\bf 156}, no.
742, 74 pp., 2002.

\bibitem{G2} K. Gr\"{o}chenig: {\it Foundations of
time-frequency analysis.} Birkh\"{a}user, Boston, 2000.

\bibitem{HW} C. Heil and D. Walnut: {\it Continuous and
discrete wavelet transforms.} SIAM Review {\bf 31}, 628--666, 1989.

\bibitem{H}
M. Holschneider: {\it Wavelets. An analysis tool.} Oxford
Mathematical Monographs. Oxford Science Publications. The Clarendon
Press, Oxford University Press, New York, 1995.

\bibitem{J} A. J. E. M. Janssen: {\it The duality condition for
Weyl-Heisenberg frames.} In ``Gabor analysis: theory and
applications'' (eds.\ H.G. Feichtinger and T. Strohmer).
Birkh\"auser, Boston, 1998.

\bibitem{L} R. S. Laugesen: \textit{Gabor dual spline windows.}
Preprint, 2008.
\\
{\small \url{www.math.uiuc.edu/~laugesen/}}

\bibitem{RoSh5} A. Ron and Z. Shen: {\it Frames and stable bases for
shift-invariant subspaces of $L^2(\mr^d)$.} Canad. J. Math. {\bf
47}, 1051--1094, 1995.

\bibitem{Wa} D. Walnut: {\it Continuity properties of the Gabor
frame operator.} J. Math. Anal.Appl. {\bf 165}, 479--504, 1992.

\bibitem{Y} R. Young: {\it An introduction to nonharmonic
Fourier series.} Academic Press, New York, 1980 (revised first
edition 2001).

\end{thebibliography}
\end{document}